\documentclass[12pt]{amsart}
\usepackage{hyperref}
\usepackage{amssymb, mathrsfs,amsmath}
\usepackage{latexsym,array}
\usepackage{amsfonts}
\usepackage{shadow}

\newtheorem{Pa}{Paper}[section]
\newtheorem{Tm}[Pa]{{\bf Theorem}}

\newtheorem{Pn}[Pa]{{\bf Proposition}}

\date{}
\author[D. Alpay]{Daniel Alpay}
\author[P. Jorgensen]{Palle Jorgensen}
\address{(DA) Department of Mathematics
\newline
Ben Gurion University of the Negev \newline P.O.B. 653,
\newline
Be'er Sheva 84105, \newline ISRAEL} \email{dany@math.bgu.ac.il}
\address{(PJ)
Department of Mathematics\newline 14 MLH \newline The University
of Iowa Iowa City,\newline IA 52242-1419 USA}
\email{jorgen@math.uiowa.edu}
\author[I. Lewkowicz]{Izchak Lewkowicz}
\address{(IL) Department of Electrical Engineering
\newline
Ben Gurion University of the Negev \newline P.O.B. 653,
\newline
Be'er Sheva 84105, \newline ISRAEL } \email{izchak@ee.bgu.ac.il}
\author[I. Marziano]{Itzik Marziano}
\address{(IM) Department of Mathematics
\newline
Ben Gurion University of the Negev \newline P.O.B. 653,
\newline
Be'er Sheva 84105, \newline ISRAEL } \email{martzian@bgu.ac.il}

\thanks{D. Alpay thanks the
Earl Katz family for endowing the chair which supported his
research. The work was done in
part while the second named author visited Department of
Mathematics, Ben Gurion University of the Negev, supported by a
BGU distinguished visiting scientist program. Support and
hospitality is much appreciated. We acknowledge discussions with
colleagues there, and in the US, Dorin Dutkay, Myung-Sin Song,
and Erin Pearse.}
\title
[Cuntz relations and interpolation] {
Representation formulas for Hardy space functions through the Cuntz relations
and new
interpolation problems}

\begin{document}
\begin{abstract}
We introduce connections between the Cuntz relations and
the Hardy space $\mathbf H_2$
of the open unit disk $\mathbb D$. We then use them to solve a
new kind of multipoint interpolation problem in $\mathbf H_2$, where for instance,
only a linear combination of the values of a function
at given points is preassigned, rather than the values at the points themselves.
\end{abstract}
\keywords{Cuntz relations, Leech's theorem, Schur analysis}
\subjclass{42C40, 47A57, 93B28}
\maketitle
\tableofcontents
\section{Introduction}
\setcounter{equation}{0}
One motivation for studying representations of the Cuntz relations
comes from signal processing, sub-band filters, and their applications to wavelets.
This falls within a larger context of multiscale problems, see for example \cite{BrJo02a}.
In this work we study the Cuntz relations in a different context,
and introduce connections between them and
the Hardy space $\mathbf H_2$ of the open unit disk $\mathbb D$.
We prove in particular
the following results: Let $b$ be a finite Blaschke product of
degree $M$, and let $e_1,\ldots, e_M$ be an orthonormal basis of
$\mathbf H_2\ominus b\mathbf H_2$. A function $f$ belongs to
$\mathbf H_2$ and has norm less or equal to $1$ if and only if it
can be written as
\begin{equation}
\label{eq:re222}
f(z)=\sum_{j=1}^Me_j(z)f_j(b(z)),
\end{equation}
where $f_1,\ldots, f_M\in\mathbf H_2$, are uniquely defined, and
are such that
\begin{equation}
\label{inequ23445}
\|f\|_{{\mathbf H}_2}^2=\sum_{j=1}^M
\|f_j\|_{{\mathbf H}_2}^2.
\end{equation}
From now on, we denote by $\|f\|_{2}$ the norm of an element of
$\mathbf H_2$.
Using  Leech's factorization theorem (see Section \ref{sec:leech} below),
we prove that, equivalently,
$f$ belongs to
$\mathbf H_2$ and has norm less or equal to $1$ if and only if it
can be written as
\begin{equation}
\label{eq:rep22}
f(z)=\frac{\sum_{j=1}^M
e_j(z)\sigma_{1j}(b(z))}{1-b(z)\sigma_2(b(z))},
\end{equation}
where
\begin{equation}
\label{eq:part11111}
\sigma=\begin{pmatrix} \sigma_{11}\\ \vdots\\ \sigma_{1M}\\
\sigma_2\end{pmatrix}
\end{equation}
is a Schur function, i.e. analytic and contractive in $\mathbb D$.\\

Representation \eqref{eq:re222} allows us to solve various
interpolation problems in $\mathbf H_2$ by translating them into
tangential interpolation problems at one point (in fact at the
origin) in $\mathbf H_2^M$. The solution of this latter problem,
or more generally, of the bitangential interpolation problem in
$\mathbf H_2^{p\times q}$ is well known. See for instance
\cite{MR96h:47020,ab6}.\\

Similarly, the representation \eqref{eq:rep22} allows us to solve
various interpolation problems in $\mathbf H_2$ by translating
them into tangential interpolation problems at {\sl one} point
(here too, in fact at the origin), for $\mathbb C^{M+1}$-valued
Schur functions, whose solution is well known. See for instance
\cite{bgr} for the general
bitangential interpolation problem for matrix-valued Schur functions.\\

We now illustrate these points. First note that, for $b$ the Blaschke
product with zeroes the points
$a_1,\ldots, a_M$,  \eqref{eq:re222} leads to
\begin{equation}
\label{lyon11}
f(a_\ell)=\sum_{j=1}^Me_j(a_\ell)f_j(0),\quad
\ell=1,\ldots M.
\end{equation}
For preassigned values of $f(a_\ell)$, $\ell=1,\ldots, M$, this
reduces the Nevanlinna-Pick interpolation problem for $M$ points
in $\mathbf H_2$ to a tangential interpolation problem at the
origin for functions in $\mathbf H_2^{M}$, whose solution, as
already mentioned, is well known. The novelty in the present
paper is by exploiting the above reduction scheme to solve
multipoint interpolation problems in $\mathbf H_2$. For example,
consider the following problem:\\

{\bf Problem:} {\sl Given $M$ points $a_1,\ldots, a_M$ in $\mathbb
D$ and $u=\begin{pmatrix}u_1&u_2&\cdots&
u_M\end{pmatrix}\in\mathbb C^{1\times M}$ and $\gamma\in\mathbb
C$, find all $f\in\mathbf H_2$ such that
\begin{equation}
\label{wertyui}
\sum_{\ell=1}^Mu_\ell f(a_\ell)=\gamma.
\end{equation}
}

{\bf Solution using \eqref{eq:re222}:}
 It follows from \eqref{lyon11} that
\begin{equation}
\label{eq:toulouse1}
\begin{split}
\sum_{\ell=1}^Mu_\ell
f(a_\ell)&=\sum_{j=1}^M(\sum_{\ell=1}^Mu_\ell e_j(a_\ell))f_{j} (0)\\
&=\sum_{j=1}^M v_jf_{j} (0),
\end{split}
\end{equation}
with
\begin{equation}
\label{vj}
v_j=\sum_{\ell=1}^Mu_\ell e_j(a_\ell),\quad j=1,\ldots, M.
\end{equation}
For preassigned value of the left side of \eqref{eq:toulouse1}
this is a classical tangential interpolation problem for $\mathbb
C^{M}$-valued functions with entries in the Hardy space. Let
$v=\begin{pmatrix} v_1&v_2&\cdots&v_M\end{pmatrix}\in\mathbb
C^{1\times M}$. Assuming $vv^*\not=0$ we have that the set of
solutions is given by
\[
\begin{pmatrix}f_1(z)\\ f_2(z)\\ \vdots \\ f_M(z)\end{pmatrix}=
\gamma\frac{v^*}{vv^*}+\left(I_M+(z-1)\frac{v^*v}{vv^*}\right)
\begin{pmatrix}g_1(z)\\ g_2(z)\\ \vdots \\ g_M(z)\end{pmatrix},
\]
where $g_1,\ldots, g_M\in\mathbf H_2$ and
\[
\sum_{\ell=1}^M\|f_j\|_2^2=\frac{|\gamma|^2}{vv^*}+\sum_{\ell=1}^M\|g_j\|_2^2.
\]
It follows from \eqref{eq:re222} that a function $f\in\mathbf H_2$ satisfies
\eqref{wertyui} if and only if it can be written as
\[
f(z)=\frac{\gamma}{vv^*}
\sum_{j=1}^Me_j(z)v_j^*+\begin{pmatrix}e_1(z)&e_2(z)&\cdots &e_M(z)\end{pmatrix}B(z)
\begin{pmatrix}g_1(z)\\ g_2(z)\\ \vdots \\ g_M(z)\end{pmatrix},
\]
where we have denoted $B(z)=\left(I_M+(z-1)\frac{v^*v}{vv^*}\right)$.
Note that $B$ is an elementary Blaschke factor,
with zero at the origin.\\

{\bf Solution using \eqref{eq:rep22}:} In the case of
representation \eqref{eq:rep22}, we have similarly
\begin{equation}
\label{lyon}
f(a_\ell)=\sum_{j=1}^Me_j(a_\ell)\sigma_{1j}(0),\quad
\ell=1,\ldots, M.
\end{equation}
For preassigned values of $f(a_\ell)$, $\ell=1,\ldots, M$, this
reduces the Nevanlinna-Pick interpolation problem for $M$ points
in $\mathbf H_2$ to a tangential interpolation problem at the
origin for matrix-valued Schur functions \eqref{eq:part11111}. As
in the previous discussion we exploit the above reduction scheme
to interpolation problem in the Schur class for multipoint
interpolation problems. For example, in the case of the
interpolation constraint \eqref{wertyui}, it follows from
\eqref{lyon} that
\begin{equation}
\label{eq:toulouse}
\begin{split}
\sum_{\ell=1}^Mu_\ell
f(a_\ell)&=\sum_{j=1}^Mu_\ell(\sum_{\ell=1}^Me_j(a_\ell))\sigma_{1j} (0)\\
&=\sum_{j=1}^M v_j\sigma_{1j} (0),
\end{split}
\end{equation}
with $v_1,\ldots, v_M$ as in \eqref{vj}.
For preassigned value of the left side of
\eqref{eq:toulouse} this is a classical tangential interpolation problem
for $\mathbb C^{M+1}$-valued Schur functions.\\

Problems of the form \eqref{lyon} have been studied for $M=2$,
under the name multipoint interpolation problem, in \cite{MR99h:47021}.
In that paper, an involution $\varphi$ of the open unit disk which maps $a_1$
into $a_2$ is used. Then one notes that the
function
\[
F(z)=\begin{pmatrix}f(z)\\ f(\varphi(z))\end{pmatrix}
\]
satisfies the symmetry
\[
F(\varphi(z))=JF(z),\quad{\rm where}\quad J=\begin{pmatrix}0&1\\ 1&0\end{pmatrix}.
\]
This reduces the interpolation problem in $\mathbf H_2$ to an
interpolation problem with symmetries in $\mathbf H^2_2$.
Unfortunately this method  does not extend to the case $M>2$. For a
related interpolation problem (for Nevanlinna functions), see also
\cite{MR2116455}, where the $n$-th composition of the
map $\varphi$ is equal to the identity map: $\varphi^{\circ n}(z)=z$.\\

In the present paper, we use a decomposition of elements in
$\mathbf H_2$ associated with isometries defined from $b$ and
which satisfy the Cuntz relation. The representation of Hardy
functions, proved in \cite{ABP}, plays a major role in the
reduction to interpolation problems in the setting of Schur
functions. To ease the notation, we set the discussion in the
framework of scalar-valued functions, but the paper itself (as
well as \cite{ABP}) is developed for matrix-valued functions.
Besides being a key player in complex analysis, the Hardy space
$\mathbf H_2$ of the open unit disk plays an important role in
signal processing and in the theory of linear dynamical systems.
An element $f$ in $\mathbf H_2$  can be described in (at least)
three different ways. In terms of: $(1)$ power series, $(2)$
integral conditions, or $(3)$, a positive definite kernel. More
precisely, in case $(1)$ one sees $f$ as the $z$-transform of a
discrete signal with finite energy, that is the $z$-transform of
a sequence $(f_n)_{n\in\mathbb N_0}\in\ell_2$:
\[
f(z)=\sum_{n=0}^\infty f_n z^n,\quad \|f\|_{\mathbf
H_2}^2\stackrel{\rm def.}{=} \sum_{n=0}^\infty|f_n|^2<\infty,
\]
In case $(2)$, one expresses the norm (in the equivalent way) as
\[
\|f\|_{2}^2=
\frac{1}{2\pi}\sup_{r\in(0,1)}\int_0^{2\pi}|f(re^{it})|^2dt<\infty,
\]
and sees $f$ as the transfer function (filter) of a
$\ell_1$-$\ell_2$ stable linear system. See \cite{MR1200235}. In
case $(3)$ we use the fact that $\mathbf H_2$ is the reproducing
kernel Hilbert space with reproducing kernel $\frac{1}{1-zw^*}$.
From the characterization of elements in a reproducing kernel
Hilbert space, a function $f$ defined in $\mathbb D$ belongs to
$\mathbf H_2$ if and only if for some $M>0$, the kernel
\begin{equation}
\label{kfh2}
\frac{1}{1-zw^*}-\frac{f(z)f(w)^*}{M}
\end{equation}
is positive definite there. The smallest such $M$ is
$\|f\|_{2}^2$. For $M=1$, rewriting \eqref{kfh2} as
\[
\frac{a(z)a(w)^*-h(z)h(w)^*}{1-zw^*}, \quad{\rm with}\quad
a(z)=\begin{pmatrix}1&-zf(z)\end{pmatrix},
\]
and using Leech's factorization theorem (see next section), it was proved in
\cite{ABP} that $f$ admits a (in general not unique)
representation of the form
\begin{equation}
\label{eq:rty}
f(z)=\frac{\sigma_1(z)}{1-z\sigma_2(z)},
\end{equation}
where $\sigma(z)=\begin{pmatrix}\sigma_1(z)\\ \sigma_2(z)\end{pmatrix}$ is
analytic and contractive in the open unit disk.\\

Let now $a\in\mathbb D$, and
\[
b_a(z)=\frac{z-a}{1-za^*}.
\]
In \cite{am_ieot} was proved that the map
\begin{equation}
T_af(z)=\frac{\sqrt{1-|a|^2}}{1-za^*}f(b_a(z))
\label{eqta}
\end{equation}
is from $\mathbf H_2$ onto itself and unitary. In the present
paper we replace $b_a$ by an arbitrary finite Blaschke product,
and define a counterpart of the operator $T_a$. If $M={\rm
deg}~b$, we now have, instead of the unitary map, $T_a$ a set of
isometries $S_1,\ldots, S_M$ in $\mathbf H_2$,  defined as
follows: Take $e_1,\ldots, e_M$ be an orthonormal basis of the
space $\mathbf H_2\ominus b \mathbf H_2$. Then,
\begin{equation}
(S_jh)(z)=e_j(z)h(b(z)),\quad h\in\mathbf H_2,
\label{eq:S_jsc}
\end{equation}
with $S_1,\ldots, S_M$ satisfying the Cuntz relations:
\begin{eqnarray}
\sum_{j=1}^M S_jS_j^*&=&I_{\mathbf H_2}, \label{cuntz11}
\\
S_j^*S_k&=&\begin{cases}I_{\mathbf H_2},\quad {\rm if}\quad j=k,\\
0,\quad\hspace{10mm} {\rm otherwise}.\end{cases}
\end{eqnarray}
It follows from these relations that every element $f\in\mathbf
H_2$ can be written in a unique way as \eqref{eq:re222}:
\begin{equation*}
f(z)=\sum_{j=1}^M e_j(z)f_j(b(z)),
\end{equation*}
where the $f_j\in\mathbf H_2$ and satisfy \eqref{inequ23445}
\[
\|f\|_2^2=\sum_{j=1}^M\|f_j\|_2^2.
\]
We note that $S_1,\ldots, S_M$ in \eqref{eq:S_jsc} form a finite
system of $M$ isometries with orthogonal ranges in $\mathbf H_2$,
with the sum of the ranges equal to all of $\mathbf H_2$. Thus
they define a representation of the Cuntz relations. This is a
special case of a result of Courtney, Muhly and Schmidt, see
\cite[Theorem 3.3]{CMS}. We send the reader to \cite{CMS} for a
survey of the relevant literature and in particular for a
discussion of the related papers \cite{MR0315454,MR1933352}. For
completeness, we provide a proof, in the matrix-valued case,
using reproducing kernel spaces techniques (see Section
\ref{sec4}). As already mentioned, one motivation for studying
representations of the
 Cuntz relations
comes from signal processing.
Our present
application of the Cuntz relations to Leech's problem from harmonic analysis
is entirely new.
The immediate relevance
to sub-band filters is
a careful selecting of the Cuntz isometries, one for each frequency sub-band, see \cite[Chapter 9]{Jor06a}.
In the case of wavelet applications,
the number $M$ is the scaling number characterizing the particular family of
wavelets under discussion.\\

We note that relations with the Cuntz relations in the indefinite inner product
case have been considered in
\cite{ajl2}, in the setting of de Branges Rovnyak spaces; see \cite{dbr1,dbr2}.
This suggests connections with interpolation in these
spaces (see \cite[Section 11]{ab6}, \cite{bbol2011}), which will be considered
elsewhere. We briefly discuss some of these aspects in Section \ref{sec7}.\\

Our paper is interdisciplinary, a mix of pure and applied,
and we are motivated by several prior developments and work by other authors.
This we discuss in section \ref{sec:leech}, \ref{sec:3} and \ref{sec4} below.
For the reader’s convenience, we mention here briefly some of these
connections. One motivation comes from earlier work \cite{MR96h:47020} on
two-sided, and tangential, interpolation for matrix functions;
see also references \cite{ab6} through \cite{ABP}, and \cite{MR2116455}.
We make additional connections also to to interpolation in
de Branges-Rovnyak spaces \cite{bbol2011}, to wavelet filters, see e.g.,
\cite{BrJo02a}, and to iterated function systems, see \cite{CMS}
by Courtney, Muhly, and Schmidt, and \cite{MR0315454} by Rochberg,
 to Hardy classes \cite{rosenblum, rr-univ} and to classical
harmonic analysis; see e.g., \cite{s1,s2,s3,s4}.\\

The  paper consists of six sections besides the introduction.
Sections 2 and 3 are of a review nature: In the second section we
discuss Leech's theorem and in the third section we discuss the
realization result of \cite{ABP}. In Section 4 we consider, in
the matrix-valued case, a set of operators which satisfy the Cuntz
relations, and were considered earlier in \cite{CMS} in the
scalar case. Section 5 is devoted to the proof of the matrix
version of \eqref{eq:rep22}. We use the representation theorem of
Section 5 to solve in Section 6 new types of multipoint
interpolation problems. Finally we outline in the last section
how some of the results  extend to the case of de Branges Rovnyak
spaces.

\section{Leech's theorem}
\setcounter{equation}{0}
\label{sec:leech}

As already mentioned in the introduction, we set the paper in the
framework of matrix-valued functions.
When the Taylor coefficients $f_n$ are $\mathbb C^{p\times
q}$-valued, one defines a $\mathbb C^{q\times q}$-valued quadratic
form by
\[
[f,f]\stackrel{\rm def.}{=}
\frac{1}{2\pi}\int_0^{2\pi}(f(e^{it}))^*f(e^{it})dt
=\sum_{n=0}^\infty f_n^*f_n.
\]
The space $\mathbf H^{p\times q}$ consists of the functions for
which ${\rm Tr}~[f,f]<\infty$ In \cite{ABP} a representation
theorem for elements $f\in\mathbf H_2^{p\times q}$ such that
$[f,f]\le I_q$ in terms of Schur functions was presented. See
Theorem \ref{tm:tams_95} below. Recall first that a $\mathbb
C^{p\times q}$-valued function $\sigma$ {\sl defined} in the open
unit disk is analytic and contractive in the open unit disk if and
only if the kernel
\begin{equation}
\label{Ksigma}
K_\sigma(z,w)=\frac{I_p-\sigma(z)\sigma(w)^*}{1-zw^*}
\end{equation}
is positive definite in the open unit disk. Such functions are
called Schur functions, denoted by $\mathscr S^{p\times q}$.
Given $\sigma\in\mathscr S^{p\times q}$ and a
$\mathbb C^{k\times p}$-valued function $A$ analytic in the open
unit disk, the kernel
\[
A(z)K_\sigma(z,w)A(w)^*=\frac{A(z)A(w)^*-B(z)B(w)^*}{1-zw^*}
\]
where $B=A\sigma$, is positive in $\mathbb D$. Leech's theorem
asserts that the converse holds: if $A$ and $B$ are respectively
$\mathbb C^{k\times p}$-valued and $\mathbb C^{k\times q}$-valued
functions {\sl defined} in $\mathbb D$ and such that the kernel
\begin{equation}
\label{eq:KAB}
\frac{A(z)A(w)^*-B(z)B(w)^*}{1-zw^*}
\end{equation}
is positive definite in $\mathbb D$, then there exists
$\sigma\in\mathscr S^{p\times q}$ such that $B=A\sigma$. Two
proofs of this theorem hold. The first assumes that $A$ and
$B$ are bounded in the open unit disk, and uses a commutant
lifting result of M. Rosenblum. See \cite{rosenblum} for
Rosenblum's result and \cite[Example 1, p. 107]{rr-univ},
\cite{add} for Leech's theorem. The other proof requires only
analyticity of $A$ and $B$ in $\mathbb D$, and uses tangential
interpolation theory for Schur functions, together with the
normal family theorem. One can extend these arguments to
functions of bounded type in $\mathbb D$, or even further weaken
these hypothesis. For completeness, we now outline a proof of
Leech's theorem for continuous functions $A$ and $B$. We first
recall the following: Let $N\in\mathbb N$ and let $w_1,\ldots,
w_N\in\mathbb D$, $\xi_1,\ldots, \xi_N\in\mathbb C^p$ and
$\eta_1,\ldots, \eta_N\in\mathbb C^q$. The tangential
Nevanlinna-Pick interpolation problem consists in finding all
Schur functions $\sigma\in\mathscr S^{p\times q}$ such that
\[
\sigma(w_j)^*\xi_j=\eta_j,\quad j=1,\ldots, N.
\]
The fact that the function $K_\sigma(z,w)$ defined by
\eqref{Ksigma} is positive definite in $\mathbb D$ implies that a
necessary condition for the tangential Nevanlinna-Pick
interpolation problem to have a solution is that the $N\times N$
Hermitian matrix $P$ (known as the Pick matrix) with $\ell j $ entry
\begin{equation}
\label{plj}
P_{\ell, j}=
\frac{\xi_\ell^*\xi_j-\eta_\ell^*\eta_j}{1-w_\ell w_j^*}
\end{equation}
is non negative. This condition is in fact also sufficient, and
there are various methods to describe all solutions in
terms of a linear fractional transformation. See for instance
\cite{bgr, MR2735305,MR90g:47003,MR2019351}. With this result at
hand we can outline a proof of Leech's theorem as follows: We
assume given two functions $A$ and $B$, respectively $\mathbb
C^{k\times q}$-valued and $\mathbb C^{k\times q}$-valued, continuous in
$\mathbb D$ and such that the kernel
\[
\frac{A(z)A(w)^*-B(z)B(w)^*}{1-zw^*}
\]
is positive definite there. Consider $w_1,w_2\ldots$ a countable
set of points dense in the open unit disk. The Hermitian block
matrix with $\ell j$ entry
\[
\frac{A(w_\ell)A(w_j)^*-B(w_\ell)B(w_j)^*}{1-w_\ell w^*_j},\quad
\ell,j=1,\ldots, N,
\]
is non negative, and therefore, by the above mentioned result on
Nevanlinna-Pick interpolation, there exists a Schur function
$\sigma_N\in\mathscr S^{p\times q}$ such that
\[
A(w_\ell)\sigma_N(w_\ell)=B(w_\ell),\quad \ell=1,\ldots N.
\]
To conclude the proof, one uses the normal family theorem to find
a function $\sigma\in\mathscr S^{p\times q}$ such that
$A(w_\ell)\sigma(w_\ell)=B(w_\ell)$ for $\ell\in\mathbb N$. By
continuity, this equality extends then to all of $\mathbb D$.

\section{A representation of $\mathbf H_2^{p\times 2}$ functions}
\label{sec:3}
\setcounter{equation}{0}
Leech's theorem can be used to find a representation of elements of
$\mathbf H_2$ in terms of Schur functions, as we now recall. See
\cite{ABP,abl3}. The following result is proved in \cite{ABP}. In
the scalar case, it was proved earlier by Sarason using different
methods. See \cite[p. 50]{s4}, \cite{s1,s2,s3}. In the discussion we shall
find it convenient to partition $\sigma\in
\mathscr S^{(p+q)\times q}$ as
\begin{equation}
\label{eq:sigma111}
\sigma=\begin{pmatrix}\sigma_1\\ \sigma_2\end{pmatrix},
\end{equation}
with $\sigma_1$ being $\mathbb C^{p\times q}$-valued and $\sigma_2$ being
$\mathbb C^{q\times q}$-valued.

\begin{Tm}
Let $H\in\mathbf H_2^{p\times q}$. Then the following are
equivalent:\\
$(1)$ It holds that:
\begin{equation}
\label{eq:ineq} [H,H]\le I_q,
\end{equation}
$(2)$ The kernel
\begin{equation}
\label{eq:KH} \frac{I_p}{1-zw^*}-H(z)H(w)^*
\end{equation}
is positive definite in $\mathbb D$.\\
$(3)$ There is a Schur function $\sigma\in\mathscr
S^{(p+q)\times q}$ (see \eqref{eq:sigma111}) so that
\begin{equation}
H(z)=\sigma_1(z)(I_q-z\sigma_2(z))^{-1}.
\end{equation}
\label{tm:tams_95}
\end{Tm}

The key to proof of this theorem is to note that the kernel \eqref{eq:KH}
can be rewritten in the form \eqref{eq:KAB} with
\[
A(z)=\begin{pmatrix}I_p&zH(z)\end{pmatrix}\quad{\rm and}\quad
B(z)=H(z),
\]
and apply Leech's theorem: There exists $\sigma\in\mathscr S^{(p+q)\times q}$ as in \eqref{eq:sigma111}
such that
\mbox{$A(z)\sigma(z)=B(z)$}, that is,
\[
\sigma_1(z)+zH(z)\sigma_2(z)=H(z).
\]
Equation \eqref{tm:tams_95} follows.\\

We note that an extension of the previous theorem to elements in the
Arveson space was given in
\cite[Theorem 10.3, p. 182]{abk}.

\section{The Cuntz relations
 in $\mathbf H_2^p$}
\label{sec4}
\setcounter{equation}{0}
Let $b$ be a finite Blaschke product of degree $M$, and let
\[
\mathscr H(b)=\mathbf H_2\ominus b\mathbf H_2.
\]
It is well known that this space is
finite dimensional and $R_0$-invariant, where
\[
R_0f(z)=\frac{f(z)-f(0)}{z},\quad f\in \mathscr H(b).
\]
Let
\[
\begin{pmatrix}f_1(z)&\ldots& f_M(z)\end{pmatrix}=C(I_M-zA)^{-1}
\]
denote a basis of $\mathscr H(b)$, where $(C,A)\in\mathbb
C^{1\times M}\times \mathbb C^{M\times M}$ is an observable pair, namely:
\[
\cap_{n=0}^\infty \ker CA^n=\left\{0\right\}.
\]
Since the spectrum of $A$ is inside the open unit disk, the series
\[
P=\sum_{\ell=0}^\infty A^{\ell*}C^*CA^\ell
\]
converges and $P>0$. The matrix $P$ is the Gram matrix (observability Gramian in control terminology) of the
basis $f_1,\ldots, f_M$, and satisfies:
\[
P=\frac{1}{2\pi}\int_0^{2\pi}
\begin{pmatrix}f_1(e^{it})&\ldots& f_M(e^{it})\end{pmatrix}^*
\begin{pmatrix}f_1(e^{it})&\ldots& f_M(e^{it})\end{pmatrix}dt.
\]
This matrix turns to be identical to the Pick matrix defined in \eqref{plj}.
We denote by $\mathbf H_2^p$ the Hilbert space of $\mathbb
C^p$-valued functions with entries in $\mathbf H_2$, and with
norm
\[
\|\begin{pmatrix}h_1\\ \vdots \\h_p\end{pmatrix}\|_{\mathbf
H_2^p}^2=\sum_{n=1}^p\|h_n\|^2_{\mathbf H_2}.
\]
We note that $\mathbf H_2^p$ is the reproducing kernel Hilbert
space with reproducing kernel $\frac{I_p}{1-zw^*}$.\\

We now introduce the Cuntz relations into this framework. This was treated earlier
in \cite{CMS} using different methods.

\begin{Tm}
Let $e_1,\ldots, e_M$ be an orthonormal basis of $\mathscr H(b)$
and, for $j=1,\ldots M$
\begin{equation}
(S_jh)(z)=e_j(z)h(b(z)),\quad h\in\mathbf H_2^p. \label{eq:S_j}
\end{equation}
Then the $S_j$ satisfy the Cuntz relations:
\begin{eqnarray}
\sum_{j=1}^M S_jS_j^*&=&I_{\mathbf H_2^{p}}, \label{cuntz1}
\\
S_j^*S_k&=&\begin{cases}I_{\mathbf H_2^{p}},\quad {if}\quad j=k,\\
0,\quad\hspace{10mm} {otherwise}.\end{cases} \label{cuntz2}
\end{eqnarray}
\label{tm41}
\end{Tm}

{\bf Proof:} We proceed in a number of steps. The proof of the
Cuntz identities \eqref{cuntz1}-\eqref{cuntz2} is given in STEPS 4 and 5, respectively.\\

STEP 1: {\sl The set $\mathbf H_2^p(b)$ of functions of the form
\[
F(z)=f(b(z)),\quad f\in\mathbf H_2^p,
\]
with norm
\[
\|F\|_{\mathbf H_2^p(b)}=\|f\|_{\mathbf H_2^p}.
\]
is the reproducing kernel Hilbert space with reproducing kernel
$\frac{I_p}{1-b(z)b(w)^*}$.}\\

This can be checked directly, but is also a special case of
\cite[Theorem 3.1, p. 109]{MR1197502}.\\

STEP 2: {\sl The operator $M_{e_j}$ of multiplication by $e_j$ is
an isometry from $\mathbf H_2^p(b)$ into $\mathbf H_2^p$.
Furthermore, the range of
$M_{e_j}$ and $M_{e_k}$ are orthogonal for $j\not =k$.}\\

Indeed, let $u,v\in\mathbb C^p$. It holds that
\[ \langle e_j b^nu, e_k b^mv
\rangle_{\mathbf H_2^p}=\begin{cases}v^*u\quad j=k\quad{\rm
and}\quad m=n,
\\0,\quad{\rm otherwise}.
\end{cases}
\]
We use that multiplication by $b$ is an isometry from $\mathbf
H_2$ into itself. If $n=m$ and $j=k$, the claim is clear. If
$n=m$ and $j\not=k$, this is just the orthogonality of $e_j$ and
$e_k$. If $n>m$, we have
\[
\langle e_j b^nu, e_kb^mv\rangle_{\mathbf H_2^p}= \langle e_j
b^{n-m}u,e_kv\rangle_{\mathbf H_2^p}=0,
\]
since $e_jb^{n-m}u\in b\mathbf H_2^p$ is orthogonal to $e_kv$
whose components belong to $\mathscr H(b)=\mathbf H_2\ominus b
\mathbf H_2$. The case $n<m$ is obtained by
interchanging the role of $j$ and $k$.\\

STEP 3: {\sl Let $e_1,\ldots, e_M$ denote an orthonormal basis of
$\mathscr H(b)$. Then,
\begin{equation}
\label{eq:decomp2} \mathbf H_2^p=\oplus_{j=1}^M e_j\mathbf
H_2^p(b).
\end{equation}
}

Indeed, the reproducing kernel is written in terms of the
orthonormal basis as (see for instance \cite[(6) p. 346]{aron},
\cite{Saitoh})
\begin{equation}
\label{eq:kb}
K_b(z,w)=\sum_{j=1}^Me_j(z)(e_j(w))^*.
\end{equation}
Thus
\begin{equation}
\label{eq:decomp}
\begin{split}
\frac{I_p}{1-zw^*}=\frac{I_p}{1-b(z)b(w)^*}\frac{1-b(z)b(w)^*}{1-zw^*}
=\sum_{j=1}^M k_j(z,w),
\end{split}
\end{equation}
with
\[
k_j(z,w)= \frac{e_j(z)(e_j(w))^*I_p}{1-b(z)b(w)^*}.
\]
Equality \eqref{eq:decomp} expresses the positive definite kernel
$\frac{I_p}{1-zw^*}$ as a sum of positive definite kernels. The
reproducing kernel space associated to $k_j$ is $e_j\mathbf
H_2(b)$. Therefore, \eqref{eq:decomp2} holds as a sum of vector
spaces; see \cite[p. 352]{aron}. Since, by STEP 2, $e_j\mathbf
H_2^p(b)$ is isometrically included into $\mathbf H_2^p$, the sum
is
orthogonal.\\

STEP 4: {\sl $S_j$ and $S_k$ are isometries, with orthogonal
ranges when $j\not =k$.}\\

The fact that $S_j$ is an isometry follows from STEPS 1 and 2.
Indeed the range of $S_j$ is in $\mathbf H_2^p$ by STEP 2 and
\[
\begin{split}
\|S_jh\|_{\mathbf H_2^p}^2&=\|h(b)\|_{\mathbf H_2^p(b)}^2\quad
(\mbox{\rm by STEP 2})\\
&=\|h\|_{\mathbf H_2^p}^2\quad (\mbox{\rm by STEP 1}).
\end{split}
\]
Furthermore, for $f,g\in\mathbf H_2^p$ and $j\not=k$,
\[
\langle S_j f, S_kg\rangle_{\mathbf H_2^p}=\langle M_{e_j} f(b),
M_{e_k}g(b)\rangle_{\mathbf H_2^p}=0
\]
by STEP 2.\\

STEP 5. {\sl It holds that $\sum_{j=1}^M S_jS_j^*=I_{\mathbf
H_2^p}$.}\\

Indeed, by the properties of multiplication and composition
operators in reproducing kernel Hilbert spaces, we have that, with
\mbox{$\rho_w(z)=1-zw^*$}, and $u\in\mathbb C^p$:
\[
(S_j^*\frac{u}{\rho_w})(z)=\frac{u}{\rho_{b(w)}(z)}(e_j(w))^*.
\]
Thus
\[
(\sum_{j=1}^MS_jS_j^*)(\frac{u}{\rho_w})(z)=\sum_{j=1}^M\frac{1}{1-b(z)b(w)^*}e_j(z)
(e_j(w))^*=\frac{u}{\rho_w(z)},
\]
where we have used \eqref{eq:kb} and \eqref{eq:decomp}. This ends
the proof since the closed linear span of the functions
$\frac{1}{\rho_w}$ is all of $\mathbf H_2$.
\mbox{}\qed\mbox{}\\

Thus we have the following decomposition  result for elements in
$\mathbf H_2$. When $M=1$, this result reduces
to the fact that the operator $T_a$ defined in \eqref{eqta}
is a unitary map from $\mathbf H_2$ into itself.

\begin{Tm}
Let $b$ be a finite Blaschke product of degree $M$, and let
$e_1,\ldots, e_M$ be an orthonormal basis of $\mathscr H(b)$.
Then, every element $H\in\mathbf H_2^{p\times q}$ can be written
in a unique way as
\begin{equation}
\label{ref_H}
H(z)=\sum_{j=1}^M e_j(z)H_j(b(z)),
\end{equation}
where the $H_j\in\mathbf H_2^{p\times q}$ and
\begin{equation}
\label{eq:opera}
[H,H]=\sum_{j=1}^M[H_j,H_j]
\end{equation}
\label{tm_H}
\end{Tm}

{\bf Proof:} We define operators $S_1,\ldots, S_M$ as in
\eqref{eq:S_j}. Let $H\in\mathbf H_2^{p\times q}$ and
$\xi\in\mathbb C^q$. It follows from the definition
\eqref{eq:S_j} of the $S_j$ and from \eqref{cuntz1} that
\[
H(z)\xi=\sum_{j=1}^M e_j(z)H_j\xi(z),
\]
where $H_j\in\mathbf H_2^{p\times q}$ is defined by
$H_j\xi=S_j^*(H\xi)$. Taking now into account \eqref{cuntz2} we
have
\[
\begin{split}
\xi^*[H,H]\xi&=\langle H\xi,H\xi\rangle_{\mathbf
H_2^{p}}\\&=\sum_{\ell,j=1}^M \langle S_\ell S_\ell^*(H\xi),
S_jS_j^*(H\xi)\rangle_{\mathbf H_2^{p}}
\\
&=\sum_{j=1}^M \langle  S_j^*(H\xi), S_j^*(H\xi)\rangle_{\mathbf
H_2^{p}}\\
&=\sum_{j=1}^M \langle  H_j\xi, H_j\xi\rangle_{\mathbf
H_2^{p}}\\
&=\xi^*\sum_{j=1}^M[H_j,H_j]\xi,
\end{split}
\]
and hence \eqref{eq:opera} holds.
\mbox{}\qed\mbox{}\\

\section{Representation of elements of $\mathbf H_2^{p\times q}$}
\setcounter{equation}{0}
We now present a generalization of Theorem \ref{tm:tams_95}. To this end,
we generalize \eqref{eq:part11111} to a partitioning of a matrix-valued functions
$\sigma\in\mathscr S^{(Mp+q)\times q}$ as
\begin{equation}
\label{eq:part2}
\sigma=\begin{pmatrix} \sigma_{11}\\ \vdots\\ \sigma_{1M}\\
\sigma_2\end{pmatrix},
\end{equation}
with $\sigma_{11},\ldots, \sigma_{1M}$ being
$\mathbb C^{p\times q}$-valued and $\sigma_2$ being
$\mathbb C^{q\times q}$-valued.

\begin{Tm}
Let $b$ be a preassigned finite Blaschke product, and let
$e_1,\ldots, e_M$ be an orthonormal basis of $\mathscr H(b)$. Let
$H\in\mathbf H_2^{p\times q}$. Then the
following are equivalent:\\
$(1)$ Condition \eqref{eq:ineq} holds: $[H,H]\le I_q$.\\
$(2)$ There exists $\sigma\in\mathscr S^{(Mp+q)\times q}$ such
that
\begin{equation}
\label{eq:rep2}
H(z)=\left(\sum_{j=1}^M
e_j(z)\sigma_{1j}(b(z))\right)(I_q-b(z)\sigma_2(b(z)))^{-1},
\end{equation}
where $\sigma\in\mathscr S^{(Mp+q)\times q}$ is as in \eqref{eq:part2}.
\end{Tm}

{\bf Proof:} By Theorem \ref{tm_H}, $H$ subject to
\eqref{eq:ineq} can be written in a unique way as \eqref{ref_H},
and it follows from \eqref{eq:opera} that
\[
\sum_{j=1}^M[H_j,H_j]\le I_q.
\]
Using Theorem \ref{tm:tams_95} with the function
\[
G=\begin{pmatrix} H_1\\ H_2\\ \vdots \\
H_M\end{pmatrix}\in\mathbf H_2^{Mp\times q},
\]
we see that there exists $\sigma\in\mathscr S^{(Mp+q)\times q}$
(see \eqref{eq:part2}) such that
\[
\begin{pmatrix}
 H_1(z)\\ H_2(z)\\ \vdots \\
H_M(z)\end{pmatrix}=\begin{pmatrix}\sigma_{11}(z)\\ \sigma_{12}(z)\\
\vdots\\ \sigma_{1M}(z)\end{pmatrix}(I_q-z\sigma_2(z))^{-1}.
\]
The result follows using \eqref{ref_H}.
\mbox{}\qed\mbox{}\\

The results in \cite{ABP} are a special
case of a family of interpolation problems with relaxed
constraints. See \cite{MR2195229,MR2249492}. We plan in a future
publication to consider these results in our new extended setting.

\section{New interpolation problems}
\setcounter{equation}{0}
We have outlined in the introduction the connections between 
multipoint interpolations and representations \eqref{eq:re222} and \eqref{eq:rep22}. 
We now add some details.
Interpolation problems whose solutions are outlined in the present section will be
considered in full details in a future publication.\\

{\bf The case of \eqref{eq:re222}:} We consider the following problem:
Find all functions $H\in\mathbf H_2^{p\times q}$
such that
\begin{equation}
\label{zxcvb}
\sum_{j=1}^M\xi_jH^{(j-1)}(a)=\gamma,
\end{equation}
for some pre-assigned matrices $\xi_1,\ldots, \xi_M\in\mathbb
C^{r\times p}$ and $\gamma\in\mathbb C^{r\times q}$. To solve
this problem we use \eqref{ref_H} with
\[
b(z)=\left(\frac{z-a}{1-za^*}\right)^M. 
\]
A basis of $\mathscr H(b)$ is
given by
\[
\frac{1}{1-za^*},
\frac{z}{(1-za^*)^2},\ldots,\frac{z^{M-1}}{(1-za^*)^M}.
\]
(see for instance
\cite{MR90g:47003}). Set
\begin{equation}
\label{E} E(z)=\begin{pmatrix}\frac{1}{1-za^*}I_p&
\frac{z}{(1-za^*)^2}I_p&\cdots&\frac{z^{M-1}}{(1-za^*)^M}I_p\end{pmatrix}.
\end{equation}
Since
\begin{equation}
\label{bbb}
b(a)=b^\prime(a)=\cdots=b^{(M-1)}(a)=0,
\end{equation}
and with
\[
\mathscr H(z)=\begin{pmatrix} H_1(b(z))\\ H_2(b(z))\\
\vdots \\ H_M(b(z))\end{pmatrix}\in\mathbf H_2^{Mp\times q},
\]
we have that
\[
\begin{split}
H(a)&=E(a)\mathscr H(0)\\
H^\prime(a)&=E^\prime(a)\mathscr H(0)\\
&\vdots\\
H^{(M-1)}(0)&=E^{(M-1)}(a)\mathscr H(0).
\end{split}
\]
Therefore, the interpolation problem \eqref{zxcvb} is equivalent to
\[
C\mathscr H(0)=\gamma,
\]
with $C\in\mathbb C^{r\times Mp}$ given by
\[
C=\sum_{j=0}^{M-1} \xi_j E^{(j)}(a).
\]
When $CC^*>0$, this in turn can be solved using \cite[Section 7]{ab6},
or directly, as
\[
\mathscr H(z)=C^*(CC^*)^{-1}\gamma+
\left(I_{Mp}+(z-1)C^*(CC^*)^{-1}C\right)\mathscr G(z),
\]
where $\mathscr G\in\mathbf H_2^{Mp\times q}$. The
formula for $H$ follows. We note that
\[
[\mathscr H,\mathscr H]=\gamma^*(CC^*)^{-1}\gamma+[\mathscr G,\mathscr G].
\]
The case where $CC^*$ is not invertible is solved
using pseudo-inverses.\\

{\bf The case of \eqref{eq:rep22}:} We here assume first that
$p=q=1$ and 
\[
b(z)=\prod_{\ell=1}^M\frac{z-a_\ell}{1-za_\ell^*},
\]
where the $a_\ell$ are distinct points in $\mathbb D$. We have now
\[
C=\begin{pmatrix}1&1&\cdots &1\end{pmatrix}\quad{\rm and}\quad
A={\rm diag}~(a_1^*,a_2^*,\ldots, a_M^*).
\]
Note that the pair $(C,A)$ is observable. Define
\begin{equation}
\label{pqwer}
P_{\ell j}=\frac{1}{1-a_ja_\ell^*},\quad \ell,j=1,\ldots, M.
\end{equation}
Namely we are in the case \eqref{plj} with the $\xi_j=1$ and the
$\eta_j=0$. In other words, $P$ is the Pick matrix obtained while
interpolating all the points $a_\ell$ to the origin. We mention
the papers \cite{szafraniec} and \cite{MR656249} for a related
discussion.
\begin{Pn}
Let $u=\begin{pmatrix}u_1&u_2&\cdots &u_M\end{pmatrix}\in\mathbb
C^{1\times M}$ and $\gamma\in\mathbb C$ be pre-assigned. Then the
following are equivalent:\\
$(1)$ It holds that
\[
\sum_{\ell=1}^M u_\ell
h(a_\ell)=\gamma\quad{and}\quad\|h\|_{\mathbf H_2}\le 1.
\]
$(2)$ $h$ is of the form
\[
C(I-zA)^{-1}P^{-1/2}\sigma_1(b(z))(1-b(z)\sigma_2(b(z))^{-1},
\]
where $\sigma=\begin{pmatrix}\sigma_1\\
\sigma_2\end{pmatrix}\in\mathscr S^{M+1}$ is such that
\begin{equation}
\label{eq:int_schur}
\begin{pmatrix}uP^{1/2}&0\end{pmatrix}\sigma(0)=\gamma,
\end{equation}
where $P$ is defined by \eqref{pqwer}.
\end{Pn}

When $b(z)=(\frac{z-a}{1-za^*})^M$ for some $M\in\mathbb N$ and
$a\in\mathbb D$, one obtains a different kind of interpolation
problem, as we now explain. Rewriting \eqref{eq:rep2} as
\begin{equation}
H(z)(I_q-b(z)\sigma_2(b(z)))=E(z)\sigma(b(z)),
\label{eqrty}
\end{equation}
where $E(z)$ is given by \eqref{E} with $p=1$, and
\[
\sigma_1(z)=\begin{pmatrix}\sigma_{11}(z)\\ \sigma_{12}(z)\\
\vdots\\ \sigma_{1M}(z)\end{pmatrix}.
\]
Differentiating \eqref{eqrty} $M-1$ times and taking into account
\eqref{bbb}, we obtain that
\begin{equation}
\begin{pmatrix}
H(a)\\ H^\prime(a)\\ \vdots\\ H^{(M-1)}(a)\end{pmatrix}=
\begin{pmatrix}E(a)&0_{p\times p}\\
E^\prime(a)&0_{p\times p}\\
&\vdots\\
E^{(M-1)}(a)&0_{p\times p}
\end{pmatrix}\begin{pmatrix}\sigma_1(0)\\
\sigma_2(0)
\end{pmatrix}.
\end{equation}
This allows to reduce \eqref{zxcvb} to a standard tangential interpolation problem
for Schur functions.

\section{The case of de Branges Rovnyak spaces}
\setcounter{equation}{0} \label{sec7}
Let $s$ be a Schur function. The kernel
$k_s(z,w)=\frac{1-s(z)s(w)^*}{1-zw^*}$ is positive definite in the
open unit disk, and the associated reproducing kerenl Hilbert
space will be denoted by $\mathscr H(s)$. Such spaces were
introduced and studied in depth by de Branges and Rovnyak in
\cite{dbr2}. When $s$ is an inner function (and in particular
when $s$ is finite Blaschke product), we have
\[
\mathscr H(s)=\mathbf H_2\ominus s\mathbf H_2.
\]
In general, $\mathscr H(s)$ is only contractively included in
$\mathbf H_2$. Let moreover $b$ be a finite Blaschke product. We
have
\begin{equation}
\label{decomp}
\begin{split}
\frac{1-s(b(z))s(b(w))^*}{1-zw^*}&=\frac{1-s(b(z))s(b(w))^*}{1-b(z)b(w)^*}
\frac{1-b(z)b(w)^*}{1-zw^*}\\
&=\sum_{j=1}^Me_j(z)e_j(w)^*\frac{1-s(b(z))s(b(w))^*}{1-b(z)b(w)^*},
\end{split}
\end{equation}
with $e_1,\ldots, e_M$ an orthogonal basis of $\mathscr H(b)$.
This decomposition allows us to define the operators $S_1,\ldots
, S_M$ as in \eqref{eq:S_j}, so that the following holds:

\begin{Tm}
The operators $S_1,\ldots, S_M$ are continuous from $\mathscr
H(s)$ into $\mathscr H(s(b))$ and satisfy the Cuntz relations:

\begin{eqnarray}
\sum_{j=1}^M S_jS_j^*&=&I_{\mathscr H(s(b))},
\label{dbrty}
\\
S_j^*S_k&=&\begin{cases}I_{\mathscr H(s)},\quad {\rm if}\quad j=k,\\
0,\quad\hspace{10mm} {\rm otherwise}.\end{cases}
\label{cuntzdb}
\end{eqnarray}
\end{Tm}

{\bf Proof:} We proceed in a number of steps.\\

STEP 1: {\sl The reproducing kernel Hilbert space $\mathscr
M(s,b)$ with reproducing kernel
$\frac{1-s(b(z))s(b(w))^*}{1-b(z)b(w)^*}$ consists of the
functions of the form \mbox{$F(z)=f(b(z))$}, with $f\in\mathscr
H(s)$ and norm
\[
\|F\|_{\mathscr M(s,b)}=\|f\|_{{\mathscr H(s)}}.
\]}

This follows from a direct computation.\\

STEP 2: {\sl The formula
\[
\left(T_j(k_{s(b)}(\cdot, w))\right)(z)=k_s(z,b(w))e_j(w)^*,\quad
w\in\mathbb D.
\]
defines a bounded  densely defined operator, which has an
extension to all of $\mathscr H(s(b))$, and whose adjoint is
$S_j$.}\\

This follows from the decomposition \eqref{decomp}.\\

STEP 3: {\sl \eqref{dbrty} holds.}\\

Indeed
\[
\begin{split}
\left(\left(\sum_{j=1}^M S_jS_j^*\right)k_{s(b)}(\cdot,
w)\right)(z)&=\sum_{j=1}^Me_j(z)k_s(b(z),b(w))e_j(w)^*\\
&=k_{s( b)}(z, w),
\end{split}
\]
and hence,  by continuity, equality \eqref{dbrty} holds in
$\mathscr H(s( b))$.\\

STEP 4: {\sl For $j\not =k$ we have
\begin{equation}
\label{inter}
e_j\mathscr M(s,b)\cap e_k\mathscr
M(s,b)=\left\{0\right\}
\end{equation}}

Indeed, let $\mathbf H_2(b)$ be as in STEP 1 in the proof of
Theorem \ref{tm41}. We have $\mathscr H(s)\subset \mathbf H_2$ and
hence
\[
\mathscr M(s,b)\subset\mathbf H_2(b).
\]
Thus \eqref{inter} follows from STEP 2 of that same theorem.\\

STEP 5: {\sl Let $M_{e_j}$ denote the operator of multiplication
by $e_j$. It holds that}
\begin{equation}
\label{decompqwery}
\mathscr H(s(b))=\oplus_{j=1}^M M_{e_j}
\mathscr M(s,b).
\end{equation}
This follows from the decomposition \eqref{decomp}, which implies
that the sum
\[
\mathscr H(s(b))=\sum_{j=1}^M M_{e_j} \mathscr M(s,b)
\]
holds, and from STEP 4, which insures that the sum is direct.\\

STEP 5: {\sl \eqref{cuntzdb} hold.}\\

Indeed, from STEP 4, the range of the operators $M_{e_j}$ and
$M_{e_k}$ are orthogonal for $j\not =k$, and $M_{e_j}$ is an
isometry.
\mbox{}\qed\mbox{}\\

Finally, we remark that \eqref{dbrty} leads to decompositions of
elements of the space $\mathscr H(s(b))$ in terms of elements of
the space $\mathscr H(s)$ similar to \eqref{eq:re222}: Every
element $f\in\mathscr H(s(b))$ can be written in a unique way as
\[
f(z)=\sum_{j=1}^Me_j(z)f_j(b(z)),
\]
where $f_1,\ldots, f_M\in\mathscr H(s)$. Furthermore
\[
\|f\|^2_{\mathscr H(s(b))}=\sum_{j=1}^M \|f_j\|_{\mathscr H(s)}^2.
\]
Multipoint interpolation problems can be also considered in this
setting, building in particular on the recent work of Ball,
Bolotnikov and ter Horst
\cite{bbol2011} on interpolation in de Branges Rovnyak spaces. 
This will be developed in a separate publication.\\

Our paper is meant as an interdisciplinary contribution, and
it involves an approach to filters and to operators having its
genesis in many different fields, both within mathematics and within
engineering. We hope that we have succeeded at least partially in communicating
across traditional lines of division separating these fields. As a result our
listed references included below is likely to be incomplete. We thank in particular Professor Paul Muhly
for improving our reference list.\\

With apologies to Goethe and to Frenchmen:\\

 {\sl Mathematicians are like Frenchmen: whatever you say to them they
translate into their own language and forthwith it is something entirely different}\\

Johann Wolfgang von Goethe.

\bibliographystyle{plain}
%
\def\cprime{$'$} \def\lfhook#1{\setbox0=\hbox{#1}{\ooalign{\hidewidth
  \lower1.5ex\hbox{'}\hidewidth\crcr\unhbox0}}} \def\cprime{$'$}
  \def\cprime{$'$} \def\cprime{$'$} \def\cprime{$'$} \def\cprime{$'$}

\end{document}